\RequirePackage{fix-cm}
\documentclass[smallextended]{svjour3}       
\smartqed  
\usepackage{graphicx}

\usepackage{epsfig,amsmath,amssymb
}

\usepackage{bm}

\usepackage{pgf}
\usepackage{tikz}
\usetikzlibrary{arrows,automata}
\usepackage[latin1]{inputenc}
\usepackage{verbatim}
\usepackage{pgfplots}

\usepackage{kpfonts}
\usepackage{stackengine}
\usepackage{calc}
\newlength\shlength
\newcommand\xshlongvec[2][0]{\setlength\shlength{#1pt}%
  \stackengine{-5.6pt}{$#2$}{\smash{$\kern\shlength%
    \stackengine{7.55pt}{$\mathchar"017E$}%
      {\rule{\widthof{$#2$}}{.57pt}\kern.4pt}{O}{r}{F}{F}{L}\kern-\shlength$}}%
      {O}{c}{F}{T}{S}}

\newcommand{\be}{\begin{equation}}
\newcommand{\ee}{\end{equation}}

\newcounter{unnumber}

\newtheorem{thrm}{Theorem}

\newtheorem{cor}{Corollary}

\newtheorem{prf}[unnumber]{Proof}

\begin{document}

\title{Elementary functions solutions to the Bachelier model generated by Lie point symmetries
}


\author{Evangelos Melas
}


\institute{E. Melas \at
              University of Athens \\
              Department of Economics \\
              Unit of Mathematics and Informatics \\
              Sofokleous 1 \\
              10559 Athens \\
              Greece \\
              Tel.: +30-210-3689403\\
              \email{emelas@econ.uoa.gr}           
}

\date{Received: date / Accepted: date}

\maketitle

\begin{abstract}
Under the recent negative interest rate situation, the Bachelier model has been attracting
attention and adopted for evaluating the price of interest rate options.
In this paper we find the Lie point symmetries of the Bachelier partial differential equation
(PDE) and use them  in order to generate new classes of denumerably infinite elementary function
solutions to the Bachelier model from elementary function solutions to it which we derived in a previous
publication.
\end{abstract}

\section{Introduction}

\label{intro}

\indent

Paul Samuelson, the first American to win the Nobel Prize in Economics, introduced Bachelier to modern financial economists in his 1972 article ``Mathematics of Speculative Price," which appeared in Mathematical Topics in Economic Theory and Computation. Samuelson said that Bachelier, ``. . . seems to have had something of a one$-$track mind. But what a track!" Samuelson attributed to Bachelier the discovery of ``Brownian motion" five years prior to Einstein's famous and much lauded similar discovery of the same phenomenon. Samuelson's encomium, most of which appears in a footnote in the article, reveals that Bachelier's work contributed to physics and mathematics before it was picked up by economists \cite{Me}.

Louis Bachelier pioneered an option pricing model in his Ph.D. thesis \cite{Ba}, marking the
birth of mathematical finance. He offered the first analysis of the mathematical properties of Brownian
motion (BM) to model the stochastic change in stock prices, and this preceded the work of Einstein  \cite{Ei}
 by five years. His analysis also precursors what is now known as the efficient market hypothesis \cite{Sc}.
See \cite{Su} for Bachelier's contribution
to financial economics and \cite{Co} for a review of his life and achievements.

Owing to the celebrated Black$-$Scholes$-$Merton (BSM)   model \cite{BSM1}, \cite{BSM2}, \cite{BSM3}, \cite{BSM4}   and the
fact that the arithmetic BM allows negative asset prices, the Bachelier model has been forgotten as a part
of history until recently. Ironically, the model gained attention again in the twenty$-$first century because
it can deal with negative asset prices, which was considered its limitation. The negative interest rates
observed in some developed countries after the 2008 global financial crisis forced fixed$-$income trading
desks to reconsider their option pricing models. The spread of COVID$-$19 led to lockdowns worldwide
and an extremely sharp drop in the global demand for oil. Consequently, in April 2020, the price of
oil futures contracts became sharply negative for the first time in history. In response, the Chicago
Mercantile Exchange (CME) and Intercontinental Exchange (ICE) changed their models for oil futures
options from the BSM to the Bachelier model to handle the negative prices \cite{CME},  \cite{ICE} .

In fact, the attention on the Bachelier model dates back to the pre$-$2008 crisis era, even when the
fear of negative prices was negligible. The Bachelier model has been widely used at least in the fixed
income markets$-$swaptions are quoted and risk$-$managed by Bachelier volatility.
Other than in fixed income markets, the Bachelier model was naturally adopted when the
underlying price can assume negative values. For example, it has been used for spread options (i.e., the
option on the price difference of two assets) \cite{Po} and year$-$on$-$year inflation \cite{Ke}.

In \cite{MeL} by studying the differential Galois group of  Bachelier ODE we obtained
four classes of denumerably infinite elementary function solutions to the Bachelier model which are expressed
via products of polynomials and exponential functions.
In this paper we find the Lie point symmetries of the Bachelier PDE
and we use these symmetries in order to generate new classes of solutions from the solutions
we derived in \cite{MeL}. We anticipate that both the solutions we found in \cite{MeL}
and the solutions we find in this paper describe various types of new interesting financial instruments.

This paper is organised as follows: In Section \ref{model} we give the bare essentials of the  Bachelier model and we derive the
Bachelier PDE. In Section \ref{Lie} we find the Lie point symmetries of the Bachelier PDE.
In Section \ref{Elementary} we revise the elementary functions solutions to the Bachelier model we found in \cite{MeL}.
In Section \ref{New} we generate  new classes of elementary function solutions  to the  Bachelier model.
In Section \ref{future} we comment on our results and we outline prospects for future research.







\section{The Bachelier model and the Bachelier PDE}

\label{model}

\indent

Following closely \cite{Te}, we give the bare essentials of the Bachelier model and the associated Bachelier PDE.

\subsection{The Bachelier model}

\indent

We consider an economy in which there is a stock without dividend payments, a zero
coupon bond (hereinafter called ``bond") and an option whose underlying asset is this
stock.
The price of stock is assumed to follow the stochastic
differential equation (SDE)
\be
dS_{t}=\mu S_{t} dt + \sigma d W_{t}^{\mathcal P},
\ee
where $\mu$ and $\sigma$ are constants and $W_{t}^{\mathcal P}$ is a
Brownian motion on a filtered probability space ($\Omega, \mathcal F, \mathcal P; \mathcal F_{t}$).
That is, $S_{t}$ follows Ornstein$-$Uhlenbeck process.

The bond price is modeled by
\be
\frac{dP(t,T)}{P(t,T)}=rdt, \ P(T,T)=1,
\ee
where $P(t,T)$  represents the value of a bond at $t$ with maturity $T$ and $r$ is a continuous compound interest rate,
which is assumed to be constant. At this time, when we choose
$P(t,T)$ as the numeraire, the SDE of $S_{t}$ under the forward measure $Q_{T}$ is
\be
dS_{t}=r S_{t} dt + \sigma d W_{t}^{Q_{T}}.
\ee

We consider a European call option on the stock with maturity date $T$ and strike
price $K$ . Thus, a payoff at $T$ is
\be
max(S_{T}-K)\equiv (S_{T}-K)^{+}.
\ee
There are no payments from contract point to maturity other than option premiums
and maturity pay$-$offs. In addition, there are no transaction costs and no restrictions on
short$-$selling, and investors can continuously restructure their portfolio.

The Bachelier PDE is derived from the Bachelier model in subsection \ref{Bachelier}
following an argument similar to the one followed by
Black $\&$ Scholes \cite{BSM1}, \cite{BSM2} and Merton \cite{BSM3} (see also \cite{Br}).

\subsection{The Bachelier PDE}
\label{Bachelier}

\indent

Let $Y_{t}$ denote the price of a European call
option with strike price $K$ and maturity $T$. We assume that it can be written as a twice
continuously differentiable function of stock prices $S_{t}$ and time $t$. That is,
\be
Y_{t}=C(S_{t},t).
\ee
Applying Ito formula to $Y_{t},$ we obtain
\be
\label{eq1}
dY_{t}=\left ( r S_{t} \frac{\partial C}{\partial S}
+ \frac{\partial C}{\partial t}    + \frac{1}{2} \sigma^{2}
\frac{\partial^{2} C}{\partial S^{2}}\right )dt +
\frac{\partial C}{\partial S} \sigma d W_{t}^{Q_{T}}.
\ee

Let $a_{t}$ denote the number of stocks held at time $t$, and let $b_{t}$ denote the number of bonds held at time $t$.
We
assume that there is a self$-$financing trading strategy $a_{t},b_{t}$  that satisfies
\be
\label{eq}
a_{t} S_{t} + b_{t} P(t,T)=Y_{t}, \ \forall t \in \left[0,T  \right ].
\ee
By linearity of the stochastic integral and the self$-$financing condition, we obtain
\begin{eqnarray}
dY_{t}  &  =   &  a_{t} d S_{t} +   b_{t} dP(t,T) \nonumber \\
\label{eq2}
&  =   & \left ( a_{t} r S_{t} +  b_{t} r  P(t,T)     \right )dt +
\alpha_{t} \sigma d W_{t}^{Q_{T}}.
\end{eqnarray}

Therefore, two expressions are obtained for $dY_{t}$, and the drift and diffusion terms of (\ref{eq1})
and (\ref{eq2}) must be equal. First, comparing the diffusion terms, we obtain
\be
a_{t}= \frac{\partial C}{\partial S}(S_{t},t).
\ee
On the other hand, from (\ref{eq}) we get
\be
b_{t}=\frac{1}{P(t,T)} \left( C(S_{t},t) - S_{t} \frac{\partial C}{\partial S}(S_{t},t)     \right ).
\ee

Consequently, comparison of the drift terms in (\ref{eq1})  and (\ref{eq2})  yields the following PDE
\be
\label{pde}
r S  \frac{\partial C}{\partial S} +
\frac{1}{2} \sigma^{2}
\frac{\partial^{2} C}{\partial S^{2}}  +
\frac{\partial C}{\partial t}  - r C=0.
\ee
This is the Bachelier PDE.

\section{Lie symmetry analysis}

\label{Lie}

\indent


Let
\be
F=r S  \frac{\partial C}{\partial S} +
\frac{1}{2} \sigma^{2}
\frac{\partial^{2} C}{\partial S^{2}}  +
\frac{\partial C}{\partial t}  - r C.
\ee

According to the Lie symmetry theory \cite{Ol}, \cite{Bl},    \cite{Ar},  \cite{Ca}, \cite{Hy} the construction of the symmetry group of (\ref{pde})
is equivalent to determination of its
infinitesimal transformation
\be
\xi= \mathcal T(t,S,C) \frac{\partial}{\partial t} +  \mathcal S(t,S,C) \frac{\partial}{\partial S}   +
\mathcal C(t,S,C)    \frac{\partial}{\partial C}.
\ee
So the second prolongation $pr^{2}\xi$ is
\be
pr^{2}\xi=\mathcal T \frac{\partial}{\partial t} +  \mathcal S \frac{\partial}{\partial S}   +
\mathcal C    \frac{\partial}{\partial C}+
 \mathcal C^{t}\frac{\partial}{\partial C_{t}}+
 \mathcal C^{S}\frac{\partial}{\partial C_{S}}+
  \mathcal C^{tt}\frac{\partial}{\partial C_{tt}}+
  \mathcal C^{tS}\frac{\partial}{\partial C_{tS}}+
   \mathcal C^{SS}\frac{\partial}{\partial C_{SS}}.
\ee

The invariance condition $pr^{2}\xi F|_{F=0}=0$ reads
\be
\label{maineq}
\left.         \left\{ r   \mathcal S    C_{S}  - r    \mathcal C +  r   \mathcal C^{S}    S +   \mathcal C^{t}  + \frac{1}{2}  \sigma^{2}   \mathcal C^{SS}
 \right\} \right |_{F=0}=0.
\ee
Expressions for the coefficients $\mathcal C^{S}$, $\mathcal C^{t} $, and $\mathcal C^{SS}$   can be found in many textbooks
(e.g. \cite{Bl}   pg. 67$-$68). Substituting $\mathcal C^{S}$, $\mathcal C^{t} $, and $\mathcal C^{SS}$ into
(\ref{maineq}), eliminating the quantities $u_{SS}$ by means of (\ref{pde}), and setting to zero all the coefficients of the
independent terms of the polynomial of $u$ and its partial derivatives, we obtain an overdetermined set of equations for the
unknown functions $\mathcal T, \mathcal S, \ \rm and \ \mathcal C $.   Solving the determinant equations, we then obtain

\begin{eqnarray}
\mathcal T&=& c_{1}  - \frac{e^{-2rt}    \left( c_{4} - c_{5}e^{4rt}    \right)}{2r },
\quad \quad \mathcal  S= \frac{1}{2} e^{-2rt}
\left (
2 c_{2}   e^{3rt}  + 2 c_{3}   e^{rt}
+ c_{4} S +  c_{5} S e^{4rt}
\right ),     \nonumber \\
\mathcal C &=&
\frac{1}{2\sigma^{2}}
\left ( e^{-2rt} \left ( \left (
\sigma^{2} \left ( - 2 c_{4}    +  c_{5}  e^{4rt} + 2 c_{6} e^{2rt} \right ) - 2rS
\left ( 2 c_{3}  e^{rt}
+ c_{4} S   \right ) \right ) C +  2 \sigma^{2} e^{2rt}  \psi(t,S)  \right )  \right ),
\end{eqnarray}
where $c_{1},c_{2},c_{3}, c_{4}, c_{5} ,  c_{6}$             are arbitrary
constants, $\psi(t,S)$    being an arbitrary solution of (\ref{pde}) satisfying
$$
r S  \frac{\partial \psi}{\partial S} +
\frac{1}{2} \sigma^{2}
\frac{\partial^{2} \psi}{\partial S^{2}}  +
\frac{\partial \psi}{\partial t}  - r \psi=0.
$$

The presence of these arbitrary constants leads to six$-$dimensional Lie algebra of symmetries.
A basis of this algebra is given by
\begin{eqnarray}
\xi_{1}&=&  \frac{\partial}{\partial t}, \ \ \ \ \ \ \ \  \xi_{2}=e^{rt}   \frac{\partial}{\partial S},   \nonumber \\
\xi_{3}&=&e^{-rt}   \frac{\partial}{\partial S} - \frac{2 e^{-rt}r S C}{\sigma^{2}} \frac{\partial}{\partial C}, \nonumber \\
\label{gen}
\xi_{4}&=& - \frac{1}{2}  \frac{e^{-2rt}}{r}  \frac{\partial}{\partial t}+
\frac{1}{2} e^{-2rt}  S  \frac{\partial}{\partial S}
- \frac{e^{-2rt} ( \sigma^{2} + r S^{2}   )C         }{\sigma^{2}}     \frac{\partial}{\partial C},  \\
\xi_{5}&=&  \frac{1}{2}        \frac{e^{2rt}}{r}  \frac{\partial}{\partial t}+
\frac{1}{2} e^{2rt}  S  \frac{\partial}{\partial S}
+\frac{1}{2} e^{2rt}  C     \frac{\partial}{\partial C}, \nonumber \\
\xi_{6}&=&  C   \frac{\partial}{\partial C}.
\nonumber
\end{eqnarray}

Let us now consider a point transformation G:
$(t,S,C) \longrightarrow (t^{*},S^{*},C^{*})    $ with $C^{*}(t^{*},S^{*})$
a solution of $ r S^{*}  \frac{\partial C^{*}}{\partial S^{*}} +
\frac{1}{2} \sigma^{2}
\frac{\partial^{2} C^{*}}{\partial S^{*2}}  +
\frac{\partial C^{*}}{\partial t^{*}}  - r C^{*}=0.    $
By using the six
generators in (\ref{gen}) we solve,
for each one of them,
the following system of ordinary differential equations
\begin{eqnarray}
\frac{dt^{*}}{d\lambda} &=& \mathcal T(t^{*},S^{*},C^{*}), \qquad
\frac{dS^{*}}{d\lambda} = \mathcal S(t^{*},S^{*},C^{*}), \qquad
\frac{dC^{*}}{d\lambda} = \mathcal S(t^{*},S^{*},C^{*}),
\end{eqnarray}
with initial condition
\begin{eqnarray}
t^{*}|_{\lambda=0}&=&t, \qquad  S^{*}|_{\lambda=0}=S,  \qquad C^{*}|_{\lambda=0}=C,
\end{eqnarray}
and we find that the corresponding
six
one$-$parameter Lie groups of
symmetry transformations of (\ref{pde}) are given respectively by
\begin{eqnarray}
\label{G1}
 \bm{G_{1}} \normalfont       &:&(t,S,C)  \longrightarrow (t^{*},S^{*},C^{*})=( G_{1}(t),  G_{1}(S),  G_{1}(C)        ) =
  \left( t+\varepsilon_{1},S,C\right ),
  \\
\label{G2}
\bm{G_{2}}&:&(t,S,C) \longrightarrow (t^{*},S^{*},C^{*})=( G_{2}(t),  G_{2}(S),  G_{2}(C)        ) =
 \left( t, S   +  \varepsilon_{2} e^{rt},C\right ),
 \\
\label{G3}
\bm{G_{3}}&:&(t,S,C)\longrightarrow (t^{*},S^{*},C^{*})=( G_{3}(t),  G_{3}(S),  G_{3}(C)        ) =
\left( t, S   +  \varepsilon_{3} e^{-rt},
e^{-\frac{ \varepsilon_{3} r  e^{-rt}
(\varepsilon_{3} e^{-rt} + 2 S   )                }{\sigma^{2}}}C\right ),
\\
\bm{G_{4}}&:&(t,S,C)\longrightarrow (t^{*},S^{*},C^{*})=( G_{4}(t),  G_{4}(S),  G_{4}(C)        )     \nonumber
\\
\label{G4}
&=&\left(
\frac{ln \left( e^{2rt} + \varepsilon_{4}\right )   }{2r}, \frac{e^{rt}S}{\sqrt{e^{2rt} + \varepsilon_{4}  }},
e^{- \frac{r\left( 2 \sigma^{2} t\left(e^{2rt} + \varepsilon_{4}\right ) -  \varepsilon_{4} S^{2}            \right )   }{\sigma^{2}\left(e^{2rt} + \varepsilon_{4}\right )}       }
\left(e^{2rt} + \varepsilon_{4}\right )  C \right ),
\\
\label{G5}
\bm{G_{5}}&:&(t,S,C)\longrightarrow
(t^{*},S^{*},C^{*})=( G_{5}(t),  G_{5}(S),  G_{5}(C)        )=
\left(
-\frac{ln \left( e^{-2rt} + \varepsilon_{5}\right )   }{2r},
\frac{e^{-rt}S}{\sqrt{e^{-2rt} + \varepsilon_{5}  }},
\frac{e^{-rt}C}{\sqrt{e^{-2rt} + \varepsilon_{5}  }}
\right ),
\\
\label{G6}
\bm{G_{6}}&:&(t,S,C)\longrightarrow (t^{*},S^{*},C^{*})=( G_{6}(t),  G_{6}(S),  G_{6}(C)        )=
\left(t,S, e^{\varepsilon_{6}} C\right ).
\end{eqnarray}

\section{Elementary function solutions to the Bachelier PDE}

\label{Elementary}

\indent

It is appropriate at this point to recall the following definition:
The function
\be
F(e,q;u)=\sum_{k=0}^{\infty} \frac{(e)_{k}}{(q)_{k}} \frac{u^{k}}{k!},
\ee
where the symbol
$(w)_{k}$, is the  Pochammer's symbol, and is defined by
\begin{equation}
(w)_{k}=w (w+1) ... (w+k-1),
\end{equation}
is called confluent hypergeometric function of the first kind or Kummer's function of the first kind.
When $e=-\rm m$, $\rm m$  being
a non-negative integer,
$F(e,q;u)$ is truncated and it reduces to a polynomial ${\rm P}_{\rm m}(u)=F(-\rm {m}$$,q;u)$  of degree $\rm m$
\be
\label{polsol}
F(-\rm {m}, \it q;u)=\rm 1- \frac{\rm m}{\it q} {\it u} + \frac{\rm m (\rm m -1)} {\it q (\it q+\rm 1)} \frac{{\it u}^{\rm 2}}{\rm 2!}+...+\frac{(-1)^{\rm m} \rm m !}{\it q(\it q+\rm 1)...(\it q+\rm m -1)} \frac{{\it u}^{\rm m}}{\rm m!}.
\ee

In \cite{MeL} we proved the following Theorem

\begin{thrm}
\label{th}
The Bachelier PDE
$$
r S  \frac{\partial C}{\partial S} +
\frac{1}{2} \sigma^{2}
\frac{\partial^{2} C}{\partial S^{2}}  +
\frac{\partial C}{\partial t}  - r C=0
$$
admits the following four classes of
elementary function solutions
\begin{eqnarray}
\label{class1}
C_{1,\rm n}(t,S) & = &S \ F \left (\frac{\rm n}{2 },\frac{3}{2};  - r \left ( \frac{ S}{\sigma } \right   )^{2}         \right ) e^{ \rm n \it  r t},  \\
\label{class2}
C_{2,\rm n}(t,S)&=&
F \left (\frac{\rm n}{2 },\frac{1}{2};
- r \left ( \frac{ S}{\sigma } \right   )^{2}
\right ) e^{(\rm n + 1) \it rt},  \\
\label{class3}
C_{3,\rm n}(t,S)&=&
e^{- r \left ( \frac{ S}{\sigma } \right   )^{2}} \
S  \ F \left (\frac{\rm n}{2 },\frac{3}{2}; r \left ( \frac{ S}{\sigma } \right   )^{2}  \right ) e^{-( \rm n -3) \it rt},
\hspace{0.7cm}\\
\label{class4}
C_{4,\rm n}(t,S)&=&e^{- r \left ( \frac{ S}{\sigma } \right   )^{2}}
\ F \left (\frac{\rm n}{2 },\frac{1}{2}; r \left ( \frac{ S}{\sigma } \right   )^{2}  \right ) e^{-( \rm n -2) \it rt},
\end{eqnarray}
where \rm n= \rm 0, \normalfont \it or \rm n \it is any multiple of \rm -2.
\end{thrm}

As pointed out in \cite{MeL},    since the Bachelier PDE is linear an immediate
cosequence  of Theorem \ref{th} is
the following
\begin{cor}
\label{cor}
The Bachelier PDE
$$
r S  \frac{\partial C}{\partial S} +
\frac{1}{2} \sigma^{2}
\frac{\partial^{2} C}{\partial S^{2}}  +
\frac{\partial C}{\partial t}  - r C=0
$$
admits the following
elementary function solution
\be
\label{gensol}
\mathcal C(t,S)=
\sum_{k=1}^{q} m_{k}
C_{1,\mu_{\it k}}(t,S) +
\sum_{k=1}^{r} z_{k}
C_{2,\zeta_{\it k}}(t,S)+
\sum_{k=1}^{t} w_{k}
C_{3,\omega_{\it k}}(t,S)+
\sum_{k=1}^{j} d_{k}
C_{4,\delta_{\it k}}(t,S),
\ee
where $q,r,t,j$ are any positive integers, equal to one or greater
than one, $m_{k},z_{k},w_{k}$ and
$d_{k}$ are arbitrary real numbers,
$\mu_{\it k},\zeta_{\it k},\omega_{\it k}$ and $\delta_{\it k}$
are
equal to zero or they are
any multiples of $-2 $,
and, $C_{\rm 1,\mu_{\it k}}(t,S),$  $C_{\rm 2, \zeta_{\it k}}(t,S),$  $C_{\rm 3, \omega_{\it k}}(t,S),$  $C_{\rm 4, \delta_{\it k}}(t,S)$
are defined respectively by (\ref{class1}), (\ref{class2}), (\ref{class3}), (\ref{class4}).

\end{cor}

\section{New solutions}

\label{New}

\indent

For convenience we write
\begin{eqnarray}
\label{class1a}
\Phi_{1,\rm n}(t,S,
C
)& = &
C -   C_{\rm 1,\rm n}(t,S), \qquad
\Phi_{2,\rm n}(t,S,C)=
C -   C_{\rm 2,\rm n}(t,S), \\
\label{class3a}
\Phi_{3,\rm n}(t,S,
C
)& = &
C
-C_{\rm 3,\rm n}(t,S), \qquad  \Phi_{4,\rm n}(t,S,C)=C-C_{\rm 4,\rm n}(t,S),
\end{eqnarray}
where $C_{\rm 1,\rm n}(t,S)$, $C_{\rm 2,\rm n}(t,S)$,  $C_{\rm 3,\rm n}(t,S)$,   $C_{\rm 4,\rm n}(t,S)$,
are defined respectively by (\ref{class1}), (\ref{class2}), (\ref{class3}), (\ref{class4}).

It is well known \cite{Ol}, \cite{Bl},    \cite{Ar},  \cite{Ca}, \cite{Hy}    that  Lie point symmetry groups of PDEs map solutions of
PDEs to different, in general, solutions of
PDEs. In the problem under consideration we have the following
\begin{thrm}
Each class of  elementary function solutions   $\Phi_{\rm q,\rm n}(t,S,C)=0$, $\rm q=\rm 1, \rm 2, \rm 3, \rm 4,$  $\rm n=\rm 0, \rm -2, \rm -4, \rm -6,...,$
to the Bachelier PDE (\ref{pde}) is mapped by the one$-$parameter  Lie group      of point symmetries  $\bm{G_{\it i}}$, $ i=\rm 1, \rm 2, \rm 3, \rm 4, \rm 5,\rm 6$,
of the Bachelier PDE (\ref{pde}), to
a class of one$-$parameter families of elementary function solutions
$\Phi_{\rm q, \rm n}  (  G_{i}(t), G_{i}(S), \newline   G_{i}(C)               )=0$,
$\rm q=\rm 1, \rm 2, \rm 3, \rm 4,$  $\rm n=\rm 0, \rm -2, \rm -4, \rm -6,...,$ $ i=\rm 1, \rm 2, \rm 3, \rm 4, \rm 5,\rm 6$,
to the Bachelier PDE (\ref{pde}).
Thus we obtain in total twenty four classes of one$-$parameter families of elementary function solutions
to the Bachelier PDE (\ref{pde}) given by
\begin{eqnarray}
\label{solution1}
\Phi_{1,\rm n}(G_{i}(t), G_{i}(S),  G_{i}(C))=0& \Leftrightarrow & G_{i}(C) -
G_{i}(S) \ F \left (\frac{\rm n}{2 },\frac{3}{2};  - r \left ( \frac{ G_{i}(S)}{\sigma } \right   )^{2}         \right ) e^{ \rm n \it  r G_{i}(t)}=0, \\
\label{solution2}
\Phi_{2,\rm n}(G_{i}(t), G_{i}(S),  G_{i}(C))=0& \Leftrightarrow & G_{i}(C) -
 \ F \left (\frac{\rm n}{2 },\frac{1}{2};  - r \left ( \frac{ G_{i}(S)}{\sigma } \right   )^{2}         \right ) e^{ \rm (n+1) \it  r G_{i}(t)}=0, \\
 \label{solution3}
\Phi_{3,\rm n}(G_{i}(t), G_{i}(S),  G_{i}(C))=0& \Leftrightarrow & G_{i}(C) -   e^{- r \left ( \frac{ G_{i}(S)}{\sigma } \right   )^{2}}
G_{i}(S) \ F \left (\frac{\rm n}{2 },\frac{3}{2};   r \left ( \frac{ G_{i}(S)}{\sigma } \right   )^{2}         \right ) e^{ -\rm (n-3) \it  r G_{i}(t)}=0, \\
\label{solution4}
\Phi_{4,\rm n}(G_{i}(t), G_{i}(S),  G_{i}(C))=0& \Leftrightarrow & G_{i}(C) -   e^{- r \left ( \frac{ G_{i}(S)}{\sigma } \right   )^{2}}
 \ F \left (\frac{\rm n}{2 },\frac{1}{2};   r \left ( \frac{ G_{i}(S)}{\sigma } \right   )^{2}         \right ) e^{ -\rm (n-2) \it  r G_{i}(t)}=0,
\end{eqnarray}
where $G_{i}(t), G_{i}(S),  G_{i}(C)$, $i=1,2,3,4,5,6,$ are given by (\ref{G1}),  (\ref{G2}),  (\ref{G3}),  (\ref{G4}),  (\ref{G5}),   (\ref{G6}).
Each class contains denumerably infinite members enumerated by $\rm n=\rm 0, \rm -2, \rm -4, \rm -6,... \ .$
Solutions (\ref{solution1}), (\ref{solution2}),  (\ref{solution3}),  (\ref{solution4}) are such that can be solved explicitly
with respect to $C.$
\end{thrm}

We note \cite{Bl}, p.85, that if for some pair ($\rm q$,  $\rm n$), $\rm q=\rm 1, \rm 2, \rm 3, \rm 4,  $
$\rm n =0,-2,-4,-6,...,$ and for some
$i=\rm 1, \rm 2, \rm 3, \rm 4, \rm 5, \rm 6,$
\be
\label{con}
\left. \xi_{i}\Phi_{\rm q,\rm n}(t,S,C)=0 \right |_{\Phi_{\rm q,\rm n}(t,S,C)=0},
\ee
then the solution $\Phi_{\rm q,\rm n}(t,S,C)=0 $
to the Bachelier PDE (\ref{pde})
is not mapped by $\bm{G_{\it i}}$,
to a  one$-$parameter family of solutions but it is mapped to itself, i.e.,
we have
\be
\Phi_{\rm q, \rm n}  (  G_{i}(t), G_{i}(S),  G_{i}(C)               )=0\overset{(\ref{con})}{\Leftrightarrow}
\Phi_{\rm q, \rm n}  (  t, S,  C               )=0.
\ee
However, this is the exception rather than the rule.

Similarly, if we write
\be
\Phi(t,S,C)=C-\mathcal C(t,S),
\ee
and,
\begin{eqnarray}
  \Omega_{1}(t,S,C)&=&C-C_{1,\mu_{\it k}}(t,S),  \qquad  \Omega_{2}(t,S,C)=C-C_{2, \zeta_{\it k}}(t,S)
\\
\Omega_{3}(t,S,C)&=&C-C_{3, \omega_{\it k}}(t,S), \qquad  \Omega_{4}(t,S,C)=C-C_{4, \delta_{\it k}}(t,S),
\end{eqnarray}
where $\mathcal C(t,S)$ is the solution of
the Bachelier PDE (\ref{pde}) defined in (\ref{gensol}),
and, $C_{\rm 1,\mu_{\it k}}(t,S),$  $C_{\rm 2, \zeta_{\it k}}(t,S),$  $C_{\rm 3, \omega_{\it k}}(t,S),$  $C_{\rm 4, \delta_{\it k}}(t,S)$,
$\mu_{\it k},\zeta_{\it k},\omega_{\it k}$ and $\delta_{\it k}$
are
equal to zero or they are
any multiples of $-2 $,
are defined respectively by (\ref{class1}), (\ref{class2}), (\ref{class3}), (\ref{class4}).

Corollary \ref{cor} leads to
\begin{cor}
\label{cor1}
The solution $\Phi(t,S,C)=0$ of the Bachelier PDE (\ref{pde}),
when at least one of the $\xi_{i}\Omega_{\rm q}$, $\rm q= \rm 1, \rm 2, \rm 3, \rm 4$  is different from zero,
is mapped by
the one$-$parameter  Lie group      of point symmetries  $\bm{G_{\it i}}$, $ i=\rm 1, \rm 2, \rm 3, \rm 4, \rm 5,\rm 6$,
of the Bachelier PDE (\ref{pde}), to
a
one$-$parameter family of elementary function solutions
\begin{eqnarray}
\Phi(G_{\it i}(t),G_{\it i}(S),G_{\it i}(C))&=&0 \Leftrightarrow G_{\it i}(C)= \mathcal C \left (  G_{\it i}(t),  G_{\it i}(S)       \right )  \Leftrightarrow     \nonumber \\
 G_{\it i}(C)&=&\sum_{k=1}^{q} m_{k}
 C_{1,\mu_{\it k}}(G_{\it i}(t),  G_{\it i}(S)) +
\sum_{k=1}^{r} z_{k}
C_{2,\zeta_{\it k}}(G_{\it i}(t),  G_{\it i}(S))  +    \nonumber \\
&&\sum_{k=1}^{t} w_{k}
C_{3,\omega_{\it k}}(G_{\it i}(t),  G_{\it i}(S))+
\sum_{k=1}^{j} d_{k}
C_{4,\delta_{\it k}}(G_{\it i}(t),  G_{\it i}(S)),
\end{eqnarray}
where $q,r,t,j$ are any positive integers, equal to one or greater
than one, $m_{k},z_{k},w_{k}$ and
$d_{k}$ are arbitrary real numbers,
$\mu_{\it k},\zeta_{\it k},\omega_{\it k}$ and $\delta_{\it k}$
are
equal to zero or they are
any multiples of $-2 $,
$C_{\rm 1,\mu_{\it k}}(t,S),$  $C_{\rm 2,\zeta_{\it k}}(t,S),$  $C_{\rm 3,\omega_{\it k}}(t,S),$  $C_{\rm 4, \delta_{\it k}}(t,S)$
are defined respectively by (\ref{class1}), (\ref{class2}), (\ref{class3}), (\ref{class4}), and,
$G_{i}(t), \ G_{i}(S),  \ G_{i}(C)$, $i=1,2,3,4,5,6,$ are given by (\ref{G1}),  (\ref{G2}),  (\ref{G3}),  (\ref{G4}),  (\ref{G5}),   (\ref{G6}).
\end{cor}

\subsection{Examples}

\subsubsection{1$^{st}$ example \rm: $\rm q=\rm 1$, $\rm n=\rm 0$, $\it i= \rm 4$.    }

\indent

We have
\be
\label{sol1}
\Phi_{1,\rm 0}(t,S,
C) = 0 \Leftrightarrow
C=   S.
\ee
(\ref{sol1}) is a solution of  (\ref{pde}). Choosing $ \it i= \rm 4 $  yields
\begin{eqnarray}
\Phi_{1,\rm 0}(G_{\rm 4}(t), G_{\rm 4}(S),  G_{\rm 4}(C) ) = 0    &\Leftrightarrow&          G_{\rm 4}(C)-   G_{\rm 4}(S)=0        \nonumber \\
&\Leftrightarrow&   e^{- \frac{r\left( 2 \sigma^{2} t\left(e^{2rt} + \varepsilon_{4}\right ) -  \varepsilon_{4} S^{2}            \right )   }{\sigma^{2}\left(e^{2rt} + \varepsilon_{4}\right )}       }
\left(e^{2rt} + \varepsilon_{4}\right )  C  - \frac{e^{rt}S}{\sqrt{e^{2rt} + \varepsilon_{4}  }}=0   \nonumber \\
\label{sol1a}
&\Leftrightarrow&
C=       \frac{    e^{ \frac{r\left( 3 \sigma^{2} t\left(e^{2rt} + \varepsilon_{4}\right ) -  \varepsilon_{4} S^{2}            \right )   }{\sigma^{2}\left(e^{2rt} + \varepsilon_{4}\right )}       }
S
}{   \left  (e^{2rt} + \varepsilon_{4}      \right  )^{\frac{3}{2}}         }.
\end{eqnarray}
(\ref{sol1a}) is a one$-$parameter family of solutions   of  (\ref{pde}). Therefore
solution (\ref{sol1})
of the Bachelier PDE (\ref{pde})  is mapped by the
one$-$parameter  Lie group      of point symmetries  $\bm{G_{\rm 4}}$, equation (\ref{G4}),
of the Bachelier PDE (\ref{pde})
to
the one$-$parameter family of solutions
(\ref{sol1a})
of the Bachelier PDE (\ref{pde}).

\subsubsection{2$^{nd}$ example \rm: $\rm q=\rm 4$, $\rm n=\rm -2$, $\it i= \rm 5$.    }

\indent

We have
\be
\label{sol2}
\Phi_{4,\rm -2}(t,S,
C) = 0 \Leftrightarrow
C=e^{4rt- r \left (\frac{s}{\sigma} \right )^{2} }
\left ( 1 - 2 r \left (\frac{s}{\sigma} \right )^{2}      \right ).
\ee
(\ref{sol2}) is a solution of  (\ref{pde}). Choosing $ \it i= \rm 5 $  yields
\begin{eqnarray}
\Phi_{4,\rm -2}(G_{\rm 5}(t), G_{\rm 5}(S),  G_{\rm 5}(C) ) = 0    &\Leftrightarrow&
G_{\rm 5}(C)-   e^{4rG_{\rm 5}(t)- r \left (\frac{G_{\rm 5}(S)}{\sigma} \right )^{2} }
\left ( 1 - 2 r \left (\frac{G_{\rm 5}(C)}{\sigma} \right )^{2}      \right )=0        \nonumber \\
\label{sol2a}
&\Leftrightarrow&
C=
\frac{
e^{ r \left (5 t - \frac{S^{2}}{\sigma^{2} \left (  1+ e^{2rt} \varepsilon_{5}   \right  )    }           \right )             }
 \sqrt{e^{-2rt}+\varepsilon_{5}}
\left ( - 2 r S^{2}+ \sigma^{2}\left (
1+ e^{2rt}\varepsilon_{5} \right) \right)
}
{\sigma^{2} \left ( 1 + e^{2rt}\varepsilon_{5}   \right )^{3}      }.
\end{eqnarray}
(\ref{sol2a}) is a one$-$parameter family of solutions   of  (\ref{pde}). Therefore
solution (\ref{sol2})
of the Bachelier PDE (\ref{pde})  is mapped by the
one$-$parameter  Lie group      of point symmetries  $\bm{G_{\rm 5}}$, equation (\ref{G5}),
of the Bachelier PDE (\ref{pde})
to
the one$-$parameter family of solutions
(\ref{sol2a})
of the Bachelier PDE (\ref{pde}).


We have
\begin{eqnarray}
\label{con1}
 \xi_{\rm 4}\Phi_{\rm 1, \rm 0}(t,S,C) &= & \left.
- \frac{  2  e^{-2rt} ( \sigma^{2} + r S^{2}   )C  +  \sigma^{2}         e^{-2rt} S       }{2    \sigma^{2}} \neq 0 \right |_{C=S},   \\
\label{con2}
 \xi_{\rm 5}\Phi_{\rm 4, \rm -2}(t,S,C) &= & \left.
\frac{1}{2 \sigma^{4}r}
\left (   e^{2rt-r \left (\frac{S}{\sigma}\right )^{2}}
\left (
14 \sigma^{2} e^{4rt}r^{2}S^{2}
- 4 e^{4rt} r^{3} S^{4} + \sigma^{4}
\left (
- 4 e^{4rt} r +
e^{r \left (\frac{S}{\sigma}\right )^{2}}
r C \right ) \right )
\right ) \neq 0
\right |_{C=e^{4rt- r \left (\frac{s}{\sigma} \right )^{2} }
\left ( 1 - 2 r \left (\frac{s}{\sigma} \right )^{2}      \right )}.
\end{eqnarray}
As a result of (\ref{con1}) and (\ref{con2}),   $\bm{G_{\rm 4}}$     and        $\bm{G_{\rm 5}}$ map
solutions
$(\ref{sol1})$  and
$(\ref{sol2})$   to
new different solutions, namely to  solutions (\ref{sol1a}) and (\ref{sol2a}) respectively.

\subsubsection{3$^{nd}$ example \rm: $q=r=t=j=2,$ $m_{\rm 1}=\rm 2, \ m_{\rm  2}=\rm 5$, $\mu_{\rm 1}=\rm 0, \ \mu_{\rm 2}=\rm -2$,
 $z_{\rm 1}=\rm 1, \ z_{\rm  2}=\rm 3$,     $\zeta_{\rm 1}=\rm 0, \ \zeta_{\rm 2}=\rm -2$,   $w_{\rm 1}=\rm 4, \it w_{\rm  2}=\rm 6$,
$\omega_{\rm 1}=\rm 0, \ \omega_{\rm  2}=\rm -2$,      $d_{\rm 1}=\rm 7, \ \it d_{\rm  2}=\rm 9$,  \ $\delta_{\rm 1}=\rm 0, \  \delta_{\rm  2}=\rm -2$.}

\indent

We have
\begin{eqnarray}
\label{sol3}
\Phi
(t,S,
C) = 0 \Leftrightarrow
C &=&   2 C_{\rm 1, \rm 0}(t,S)   +   5 C_{\rm 1, \rm -2}(t,S) +    C_{\rm 2, \rm 0}(t,S)   +   3 C_{\rm 2, \rm -2}(t,S) +   4C_{\rm 3, \rm 0}(t,S)   +   6 C_{\rm 3, \rm -2}(t,S)
\nonumber \\
&+& 7 C_{\rm 4, \rm 0}(t,S) + 9 C_{\rm 4, \rm -2}(t,S)  \nonumber \\
\Leftrightarrow C &=& 2 S + 5   S  \left (    1 +   \frac{2r}{3 }    \left (  \frac{S}{\sigma}    \right )^{2} \right )     e^{-2rt}
+ e^{rt} +  3   \left (    1 +    2 r   \left (  \frac{S}{\sigma}    \right )^{2} \right )     e^{-rt}   + 4 S e^{3rt-r  \left (   \frac{S}{\sigma}      \right )^{2} }    \nonumber \\
&+& 6 S e^{5rt-r  \left (   \frac{S}{\sigma}      \right )^{2} }
 \left (    1 -    \frac{2r}{3 }    \left (  \frac{S}{\sigma}    \right )^{2} \right )   +   7   e^{2rt-r  \left (   \frac{S}{\sigma}      \right )^{2} }
 + 9   e^{4rt-r  \left (   \frac{S}{\sigma}      \right )^{2} } \left (    1 -    2 r   \left (  \frac{S}{\sigma}    \right )^{2} \right )
\end{eqnarray}

(\ref{sol3}) is a solution of  (\ref{pde}). Choosing $ \it i= \rm 3 $  yields
\begin{eqnarray}
\Phi
(G_{\rm 3}(t), G_{\rm 3}(S),  G_{\rm 3}(C) ) = 0
\Leftrightarrow   G_{\rm 3}(C) & = &  2 G_{\rm 3}(S) + 5  G_{\rm 3}(S)   \left (    1 +   \frac{2r}{3 }    \left (  \frac{G_{\rm 3}(S)}{\sigma}    \right )^{2} \right )     e^{-2rG_{\rm 3}(t)} +   e^{r G_{\rm 3}(t)}
\nonumber \\
&+&  3   \left (    1 +    2 r   \left (  \frac{G_{\rm 3}(S)}{\sigma}    \right )^{2} \right )     e^{-rG_{\rm 3}(t)}
  + 4 G_{\rm 3}(S)    e^{3rG_{\rm 3}(t)-r  \left (   \frac{G_{\rm 3}(S)}{\sigma}      \right )^{2} }    \nonumber \\
&+& 6 G_{\rm 3}(S) e^{5rG_{\rm 3}(t)-r  \left (   \frac{G_{\rm 3}(S)}{\sigma}      \right )^{2} }    \left (    1 -    \frac{2r}{3 }    \left (  \frac{G_{\rm 3}(S)}{\sigma}    \right )^{2} \right )
 +   7   e^{2rG_{\rm 3}(t)-r  \left (   \frac{G_{\rm 3}(S)}{\sigma}      \right )^{2} } \nonumber \\
 \label{sol4}
 &+&    9   e^{4rG_{\rm 3}(t)-r  \left (   \frac{G_{\rm 3}(S)}{\sigma}      \right )^{2} } \left (    1 -    2 r   \left (  \frac{G_{\rm 3}(S)}{\sigma}    \right )^{2} \right ).
\end{eqnarray}

Equation (\ref{sol4}) gives

\begin{eqnarray}
C &=&  e^{ \frac{\varepsilon_{3} r  e^{-rt} \left (\varepsilon_{3} e^{-rt} + 2 S      \right )       }{\sigma^{2}}   }
\left (
2 \left ( S +  \varepsilon_{3}  e^{-rt}       \right ) +
\frac{ 5 e^{-5rt} \left ( \varepsilon_{3} + S e^{rt}  \right )
\left ( e^{2rt} \left ( 2 r S^{2} + 3 \sigma^{2} \right )
+ 2 \varepsilon_{3} r \left ( \varepsilon_{3} + 2 e^{rt} S
\right ) \right )
}
{3 \sigma^{2}}
\right )
\nonumber \\
&+&
e^{ \frac{\varepsilon_{3} r  e^{-rt} \left (\varepsilon_{3} e^{-rt} + 2 S      \right )       }{\sigma^{2}}   }
\left (
e^{rt} +
\frac{3 e^{-3rt}
\left (
e^{2rt} \left (
2 r S^{2} + \sigma^{2} \right )+
2 \varepsilon_{3} r
\left ( \varepsilon_{3} + 2 e^{rt} S
\right )
\right )
}
{\sigma^{2}}
\right )  \nonumber \\
&+&
e^{ \frac{\varepsilon_{3} r  e^{-rt} \left (\varepsilon_{3} e^{-rt} + 2 S      \right )       }{\sigma^{2}}   }
\left (
4 e^{2rt -  \frac{r \left (S + \varepsilon_{3} e^{-rt} \right )^{2}
}{\sigma^{2}}   }
\left (\varepsilon_{3} +  e^{rt} S  \right )
+ 6 e^{5rt -  \frac{r \left (S + \varepsilon_{3} e^{-rt} \right )^{2}
}{\sigma^{2}}   }
\left ( S +   \varepsilon_{3} e^{-rt}   \right )
\left (
1- \frac{2 r \left ( S + \varepsilon_{3}  e^{-rt}  \right )^{2}}{3 \sigma^{2}}
\right )
\right )  \nonumber \\
\label{sol3a}
&+&
e^{ \frac{\varepsilon_{3} r  e^{-rt} \left (\varepsilon_{3} e^{-rt} + 2 S      \right )       }{\sigma^{2}}   }
\left (
7 e^{2rt -  \frac{r \left (S + \varepsilon_{3} e^{-rt} \right )^{2}
}{\sigma^{2}}   }
+ 9 e^{4rt -  \frac{r \left (S + \varepsilon_{3} e^{-rt} \right )^{2}
}{\sigma^{2}}   }
\left (
1- \frac{2 r \left ( S + \varepsilon_{3}  e^{-rt}  \right )^{2}}{ \sigma^{2}}
\right )
\right ).
\end{eqnarray}

(\ref{sol3a}) is a one$-$parameter family of solutions   of  (\ref{pde}). Therefore
solution (\ref{sol3})
of the Bachelier PDE (\ref{pde})  is mapped by the
one$-$parameter  Lie group      of point symmetries  $\bm{G_{\rm 3}}$, equation (\ref{G3}),
of the Bachelier PDE (\ref{pde})
to
the one$-$parameter family of solutions
(\ref{sol3a})
of the Bachelier PDE (\ref{pde}). We note that $\xi_{\rm 3}C_{\rm q, \rm 0} \neq 0, \  $
$\rm q=\rm 1, \ \rm 2, \ \rm 3, \ \rm 4,$
on the corresponding solution surface.

\section{Conclusion
and Future Development
}
\label{future}

\indent

Analytical tractability of any financial model is an important feature. Existence of
a closed-form solution definitely helps in pricing financial instruments and calibrating
the model to market data. It also helps to verify the model assumptions, check its
asymptotic behavior and explain causality. In fact, in mathematical finance many models were
proposed, first based on their tractability, and only then by making another argument.

Our
new classes of elementary function solutions
to the Bachelier model
allow fast and accurate calculation of prices of various financial instruments
under the Bachelier process. Moreover
they will facilitate further the use of the Bachelier model as a benchmark for
the pricing of
various types of financial instruments.



The study of the financial instruments
which these elementary function solutions describe
will be useful
 for trading by using the Bachelier model in all cases.

Our new classes of elementary function solutions
to the Bachelier model have been  derived with a two$-$stage
Research Program: First, in \cite{MeL}, we obtained Liouvillian solutions
to the Bachelier model by studying its differential Galois Group. Then, in this paper,
we use the Lie point symmetries of  the Bachelier model   in order to generate new solutions
from the solutions which were derived at the first stage of the Program.

We aim to apply the same two$-$stage Research Program in order to explore the solution space
of various models in Financial Mathematics such as the Heston Model \cite{He} and its subsequent
generalizations
(the Merton$-$Garman model  \cite{BSM3},  \cite{Ga} the Chen model  \cite{Lin}, the Garch model \cite{Robert}  to name but only a few), the Longstaff$-$Schwartz model
\cite{LONG}, the Vasicek model \cite{Vas}, the Cox$-$Ingersoll$-$Ross model \cite{Cox}, the Heath$-$Platin$-$Schweizer model \cite{Heath}, e.t.c..






\end{document}